\documentclass[12pt]{article}

\usepackage{amsmath}
\usepackage{amssymb}
\usepackage{amsfonts}
\usepackage{a4}
\usepackage{graphicx}
\usepackage[abs]{overpic}
\usepackage{caption}
\usepackage{wrapfig}
\usepackage{float}
\usepackage{graphicx}
\usepackage{mathrsfs}

\newtheorem{theorem}{Theorem}[section]
\newtheorem{lemma}[theorem]{Lemma}

\newtheorem{definition}[theorem]{Definition}

\newtheorem{rem}[theorem]{Remark}

\newcommand{\Proof}{\par\noindent{\em Proof. }}
\newcommand{\eop}{\nopagebreak\hspace*{\fill}$\Box$\smallskip}
\def\diam{d}

\newcommand{\N}{\Bbb N}
\newcommand{\Z}{\Bbb Z}
\newcommand{\R}{\Bbb R}

\def\eps{\varepsilon}

\def\dist{\operatorname{dist}}

\def\XXint#1#2#3{{\setbox0=\hbox{$#1{#2#3}{\int}$}
     \vcenter{\hbox{$#2#3$}}\kern-.5\wd0}}

\usepackage{color}

\begin{document}

\begin{center}
\begin{Large}
{\bf {A Korn-type inequality in  SBD for functions with small jump sets}}
\end{Large}
\end{center}

\begin{center}
\begin{large}
Manuel Friedrich\footnote{Faculty of Mathematics, University of Vienna, 
Oskar-Morgenstern-Platz 1, 1090 Vienna, Austria. {\tt manuel.friedrich@univie.ac.at}}
\end{large}
\end{center}

\bigskip

\begin{abstract}
We present a Korn-type inequality in a planar setting for special functions of bounded deformation. We prove that for each function in $SBD^2$ with  a  sufficiently small jump set the distance of the function and its derivative from an infinitesimal rigid motion can be controlled in terms of the linearized elastic strain outside of a small exceptional set of finite perimeter.  Particularly,   the result shows that each function in $SBD^2$ has bounded variation away from an arbitrarily small part of the domain.
\end{abstract}
\bigskip

\begin{small}

\end{small}

\tableofcontents

\section{Introduction}

 The space $BD(\Omega,\R^d)$ of functions of bounded deformation, which consists of all functions $u \in L^1(\Omega,\R^d)$ whose  symmetrized distributional derivative $Eu := \frac{1}{2}((Du)^T + Du)$ is a finite $\R^{d \times d}_{\rm sym}$-valued Radon measure, has been introduced for the investigation of geometrically linear problems in plasticity theory and fracture mechanics (see \cite{Ambrosio-Coscia-Dal Maso:1997, Bellettini-Coscia-DalMaso:98, Suquet}). Variational damage or fracture problems are widely formulated in the subspace $SBD^2(\Omega,\R^d)$ (for the definition and properties of this space we refer to Section \ref{rig-sec: sub, bd} below). In the spirit of the seminal work \cite{Francfort-Marigo:1998}, the modeling essentially concentrates on the competition between elastic bulk contributions $\Vert e(u) \Vert_{L^2(\Omega)}$ given in terms of the linear elastic strain $e(u) := \frac{1}{2}(  (\nabla u)^T  + \nabla u)$ and surface terms  that assign energy contributions on the crack paths comparable to the size of the crack ${\cal H}^{d-1}(J_u \cap \Omega)$, where $J_u$ denotes the `jump set' of $u$ (see e.g. \cite{Bourdin-Francfort-Marigo:2008, Chambolle:2003, Chambolle:2004, Focardi-Iurlano:13,  SchmidtFraternaliOrtiz:2009}). 

A major additional difficulty of these problems compared to models in $SBV$ (see \cite{Ambrosio-Fusco-Pallara:2000} for the definition and basic properties of the space of special functions of bounded variation) is the lack of control on the skew symmetric part of the distributional derivative $ (Du)^T  - Du$. In fact, it is a natural and important question to analyze in which circumstances the displacement field $u$ or the absolutely continuous part of its derivative $\nabla u$ can be controlled by $\Vert e(u) \Vert_{L^2(\Omega)}$ and ${\cal H}^{d-1}(J_u)$. Apart from establishing compactness results, such properties may contribute to gain  profound  understanding of the relation between $SBD$ and $SBV$ functions which is highly desirable since in contrast to $(S)BV$ fine properties in $BD$ appear not to be well understood by now. (We refer to the recent paper \cite{Conti-Iurlano:15} for a thorough discussion  and some results in that direction.)

The key estimate providing a relation between the symmetric and the full part of the gradient is know as Korn's inequality. In its basic version,  it states that for a bounded connected Lipschitz set $\Omega$ and $p \in (1,\infty)$ there is a constant $C(\Omega,p)$ depending only on $p$ and the domain $\Omega \subset \R^d$ such that for all $u \in W^{1,p}(\Omega, \R^d)$ there is some  $A \in \R^{d \times d}_{\rm skew}$ with
\begin{align}\label{in: 1}
 \Vert \nabla u - A\Vert_{L^p(\Omega)} \le C(\Omega,p)\Vert e(u) \Vert_{L^p(\Omega)}.
 \end{align}
(See e.g. \cite{Nitsche} for a proof and \cite{Conti-Dolzmann-Muller:14, Friedrich-Schmidt:15, FrieseckeJamesMueller:02, Muller-Scardia-Zeppieri:14} for generalizations of this result into various directions.) It turns out that the statement is false in $W^{1,1}$, i.e. one can construct functions with $e(u) \in L^1(\Omega)$, but $\nabla u \notin L^1(\Omega)$ (cf. \cite{ContiFaracoMaggi:2005, KirchheimKristensen:2011, Ornstein}). On the one hand, these observations  particularly show that $BD$ is not contained in $BV$. On the other hand, it raises the natural question if in the space $SBD^2$ an estimate similar to \eqref{in: 1} can be 
established due to the higher integrability for the elastic strain $e(u)$. 

However, simple examples, e.g. in \cite{Ambrosio-Coscia-Dal Maso:1997} or the piecewise rigidity result proved in \cite{Chambolle-Giacomini-Ponsiglione:2007}, show that \eqref{in: 1} cannot hold for general functions in $SBD^2$ since the  behavior of small pieces being almost or completely detached from the bulk part of the specimen might not be controlled. In the recently appeared contributions \cite{Chambolle-Conti-Francfort:2014, Friedrich:15-1}, it has been proved that for displacement fields having small jump sets  with respect to the size of the domain the distance of the function from an infinitesimal rigid motion can be estimated in terms of the linearized elastic energy outside of a small exceptional set  $F$.  However, these  Korn-Poincar\'e-type estimates being essentially of the form 
\begin{align}\label{in: 2}
 \Vert u - (A\, \cdot + b) \Vert_{L^2(\Omega\setminus F)} \le C(\Omega)\Vert e(u) \Vert_{L^2(\Omega)}
 \end{align}
for $u \in SBD^2(\Omega,\R^d)$ and corresponding $A \in \R^{d \times d}_{\rm skew}$, $b \in \R^d$, are significantly easier as in contrast to \eqref{in: 1} no derivative is involved. The goal of the present article is to provide a generalization of \eqref{in: 2} to an inequality of Korn's-type in a planar setting  where one 
additionally controls $\nabla u$ away from a small exceptional set of finite perimeter. Our main result is the following. 

\begin{theorem}\label{th: mainkorn}
Let $\Omega\subset \R^2$ be an open, connected, bounded set with Lipschitz boundary and let $p \in [1,2), q \in [1,\infty)$. Then there is a constant $C=C(\Omega, p,q)$ such that for all $u \in SBD^2(\Omega,\R^2)$ there is a set of finite perimeter $F \subset \Omega$ with 
\begin{align}\label{eq: R2main}
{\cal H}^1(\partial^* F) \le C{\cal H}^1(J_u), \ \ \ \ |F| \le C({\cal H}^1(J_u))^2
\end{align}
and $A \in \R^{2 \times 2}_{\rm skew}$, $b \in \R^2$ such that
\begin{align}\label{eq: main estmain}
\begin{split} 
(i) & \ \ \Vert u - (A\,\cdot + b) \Vert_{L^q(\Omega \setminus F)}\le C \Vert e(u) \Vert_{L^2(\Omega)},\\
(ii) & \ \ \Vert \nabla u - A \Vert_{L^p(\Omega \setminus F)}\le C \Vert e(u) \Vert_{L^2(\Omega)},
\end{split}
\end{align}
where $e(u)$ denotes the part of the strain $Eu = \frac{1}{2}( (Du)^T  + Du)$ which is absolutely continuous with respect to  the Lebesgue measure ${\cal L}^2$,  $|F|$ stands for the ${\cal L}^2$-measure of $F$  and $\partial^* F$ is the essential boundary of $F$. 
\end{theorem}
 
 We also refer to Section \ref{rig-sec: sub, bd} for the relevant definitions.  Let us first mention that we establish the result only in two dimensions as we employ a modification technique for special functions of bounded deformation (see \cite{Friedrich:15-1}) which was only derived in a planar setting due to technical difficulties concerning the topological structure of crack geometries in higher dimensions. 

Note that  we can control the length of the boundary ${\cal H}^1(\partial^* F)$ of the exceptional set $F$ which is associated to the parts of $\Omega$ being detached from the bulk part of $\Omega$ by $J_u$. Consequently, the result is adapted for the usage of compactness theorems for $SBV$ and $SBD$ functions (see \cite{Ambrosio-Fusco-Pallara:2000, Bellettini-Coscia-DalMaso:98, DalMaso:13}).

Although the main goal of the work at hand is the derivation of the estimate for the derivative in \eqref{eq: main estmain}(ii), we also provide a generalization for the integrability exponent $q$. In \cite{Friedrich:15-1}, the exponent was restricted to $q=2$ due to the application of a $BD$ Korn-Poincar\'e inequality and in \cite{Chambolle-Conti-Francfort:2014} the arguments were based on slicing techniques similar to those used in the proof of Sobolev embeddings and led to an exponent $q= \frac{2d}{d-1}$. In the present context we obtain the estimate for $q < \infty = 2^*$ as in the usual Sobolev setting.  

To our knowledge the first inequality of this kind, i.e., also involving an estimate on $\nabla u$,
has been presented in \cite{here}, which is the preprint version of the present paper. Subsequently, the result has been extended to the critical exponent $p=2$ in \cite{Conti-Iurlano:15.2}. However, as the proof techniques are quite different, we believe that also the present article may be interesting for the community. 

As an application, we discuss that Theorem \ref{th: mainkorn} together with an approximation result shows that $SBD^2$ functions have bounded variation outside an arbitrarily small exceptional set of finite perimeter (see Theorem \ref{th: appl} below). Hereby we give another contribution to the relation between $SBV$ and $SBD$ functions which appears to go in a slightly different direction than the results presented in \cite{Conti-Iurlano:15}. We note that this statement does not immediately follow from the main theorem since a bound on $\nabla u$ does not automatically ensure that $u$ has bounded variation.

Similarly as the previously mentioned results \cite{Chambolle-Conti-Francfort:2014, Friedrich:15-1} or the $SBV$ Poincar\'e inequality \cite{DeGiorgiCarrieroLeaci:1989}, Theorem \ref{th: mainkorn} establishes an estimate only for functions whose jump set is small with respect to the size of the domain. Indeed, for larger jump sets the body may be separated into different parts of comparable size (cf. \cite{Chambolle-Giacomini-Ponsiglione:2007, Friedrich:15-2, Friedrich-Schmidt:15, FriedrichSolombrino} for related problems). In the general case, we expect a `piecewise Korn inequality' to hold,  
i.e. the body may be broken into different sets and on each connected component the distance of the displacement field from a certain infinitesimal rigid motion can be controlled. We defer the analysis of this problem, for which Theorem \ref{th: mainkorn} is a key ingredient, to a subsequent work \cite{Friedrich:15-4}. 

The paper is organized as follows. In Section \ref{rig-sec: pre} we first recall the definition and basic properties of functions of bounded variation and deformation (Section \ref{rig-sec: sub, bd}). Then in Section \ref{rig-sec: korn} we introduce the notion of \emph{John domains} being a class of sets with possibly highly irregular boundary (see e.g. \cite{John:1961, Martio-Sarvas:1978}). It is convenient to formulate Korn's and Poincar\'e's inequality for these sets since there are good criteria to obtain uniform control over the involved constants independently of the particular shape of the domain (cf. \cite{Acosta}). Finally, in Section \ref{rig-sec: modifica} we present the modification technique proved in \cite{Friedrich:15-1} which shows that after a small alteration of the displacement field and the jump set the jump heights of an $SBD$ function can be controlled solely by $\Vert e(u) \Vert_{L^2(\Omega)}$ and ${\cal H}^1(J_u)$. 
 
The rest of the paper then contains the proof of Theorem \ref{th: mainkorn}. We first establish a local estimate on a square and after a subsequent analysis of the problem near the boundary of the Lipschitz set the main theorem follows by a standard covering argument. 

In Section \ref{sec: local} we concern ourselves with the local estimate and first see that by an approximation argument (cf. \cite{Chambolle:2004}) it suffices to consider $SBD$ functions with regular jump set. The main strategy is then to modify a function with the techniques presented in Section \ref{rig-sec: modifica}. Consequently, using a Korn-Poincar\'e inequality in  $BD$ (see \cite{Kohn:82, Temam:85}), we find good approximations of the displacement field by infinitesimal rigid motions in neighborhoods of the jump set. Then drawing ideas from \cite{Friedrich-Schmidt:15} we can iteratively modify the configuration on various mesoscopic length scales to find a Sobolev function on the square which coincides with the original displacement field outside of a small exceptional set. Finally, the local estimate follows by application of the standard inequality \eqref{in: 1}.

We remark that for the nonlinear estimate \cite{Friedrich-Schmidt:15} it was not possible to gain control over  the full part of the gradient as the approximating rigid motions had to be adapted after each iteration step leading to a continual increase of the involved constant. In the present context, however, the affine mappings are found a priori and are fixed during the modification procedure whereby a bound for $\nabla u$ can be established using H\"older's and a scaled Young's inequality in the case $p<2$. Moreover, let us mention that our approach to derive the local estimate for the Korn inequality differs from the one proposed in \cite{Conti-Iurlano:15.2}, where a suitable triangulation of the domain and a corresponding linear interpolation are constructed.  

Section \ref{sec: bound} contains the main estimate at the boundary. We consider a Whitney covering of the domain and apply the result obtained in Section  \ref{sec: local}  on every square where  the length of  the jump $J_u$ is small. Hereby we can again construct a Sobolev function outside a small exceptional set $F$. For the application of \eqref{in: 1} now an additional difficulty occurs as we have to control the shape of domain. In this context, we show that choosing $F$ appropriately we find that
the complement is a John domain for a universal \emph{John constant} and therefore we can derive a uniform estimate of the form  \eqref{in: 1}. 

In Section \ref{sec: main} we then give the main proof and discuss an application to the relation of $SBV$ and $SBD$ functions. The standard examples for $(S)BD$ functions not having bounded variation are given by configurations where small balls are cut out from the bulk part with an appropriate choice of the functions on these specific sets (see e.g.  \cite{Ambrosio-Coscia-Dal Maso:1997, Conti-Iurlano:15}).  We prove that each $SBD^2$ function has bounded variation away from an arbitrarily small part of the body essentially showing that the mentioned construction provides the only way to obtain functions not lying in $SBV$.

\section{Preliminaries}\label{rig-sec: pre}

In this preparatory section  we first recall the definition and basic properties of functions of bounded variation and state Korn's and Poincar\'e's inequality for \emph{John domains}. Afterwards, we recall a result obtained in \cite{Friedrich:15-1} providing modifications of $SBD$ functions for which the jump heights can be controlled in terms of the linear elastic strain.

\subsection{Special functions of bounded variation}\label{rig-sec: sub, bd}

In this section we collect the definitions of $SBV$ and $SBD$ functions.  Let $\Omega \subset \R^d$  be an  open, bounded set  with Lipschitz boundary. Recall that the space $SBV(\Omega, \R^d)$, abbreviated as $SBV(\Omega)$ hereafter,  of \emph{special functions of bounded variation} consists of functions $u \in L^1(\Omega, \R^d)$ whose distributional derivative $Du$ is a finite Radon measure, which splits into an absolutely continuous part with density $\nabla u$ with respect to Lebesgue measure and a singular part $D^s u$. The Cantor part $D^c u$ of $D^s u$ vanishes and thus we have
$$ D^s u = [u] \otimes \xi_u {\cal H}^{d-1} \lfloor J_u, $$
where ${\cal H}^{d-1}$ denotes the $(d-1)$-dimensional Hausdorff measure, $J_u$ (the `crack path') is an ${\cal H}^{d-1}$-rectifiable set in $\Omega$, $\xi_u$ is a normal of $J_u$ and $[u] = u^+ - u^-$ (the `crack opening') with $u^{\pm}$ being the one-sided limits of $u$ at $J_u$. If in addition $\nabla u \in L^p(\Omega)$ for $1 < p < \infty$ and ${\cal H}^{d-1}(J_u) < \infty$, we write $u \in SBV^p(\Omega)$. Moreover, $SBV_{\rm loc}(\Omega)$ denotes the space of functions which belong to $SBV(\Omega')$ for every open set $\Omega' \subset \subset \Omega$.

Furthermore, we define the space $GSBV(\Omega)$ of \emph{generalized special functions of bounded variation} consisting of all ${\cal L}^d$-measurable functions $u: \Omega \to \R^d$ such that for every $\phi \in C^1(\R^d)$ with the support of $\nabla \phi$ compact, the composition $\phi \circ u $ belongs to $SBV_{\rm loc}(\Omega)$ (see \cite{DeGiorgi-Ambrosio:1988}). Likewise, we say $u \in GSBV^p(\Omega)$ for $u \in GSBV(\Omega)$ if $\nabla u \in L^p(\Omega)$ and ${\cal H}^{d-1}(J_u) < \infty$. See \cite{Ambrosio-Fusco-Pallara:2000} for the basic properties of theses function spaces. 

We now state a version of Ambrosio's compactness theorem in $GSBV$ adapted for our purposes (see e.g. \cite{Ambrosio-Fusco-Pallara:2000, DalMaso-Francfort-Toader:2005}):

\begin{theorem}\label{clea-th: compact}
Let $\Omega \subset \R^d$ open, bounded and let $1 < p < \infty$. Let $(u_k)_k$ be a sequence in $GSBV^p(\Omega)$ such that
$$\Vert \nabla u_k \Vert_{L^p(\Omega)}+ {\cal H}^{d-1}(J_{u_{k}}) + \Vert u_k \Vert_{L^1(\Omega)} \le C  $$
for some constant $C$ not depending on $k$. Then there is a subsequence (not relabeled) and a function $u \in GSBV^p(\Omega)$ such that $u_k \to u$ a.e. and $\nabla u_k \rightharpoonup \nabla u$ weakly in $L^p(\Omega)$. If in addition $\Vert u_k\Vert_{\infty} \le C$ for all $k \in \N$, we find $u \in SBV^p(\Omega) \cap L^\infty(\Omega)$.
\end{theorem}

An important subset of $SBV$ is given by the indicator functions $\chi_W$ for measurable   $W\subset \Omega$ with ${\cal H}^{d-1}( \partial^* W)< \infty$, where $\partial^* W$ denotes  the \emph{essential boundary} of $W$  (also called \emph{measure-theoretic boundary},  see \cite[Definition 3.60]{Ambrosio-Fusco-Pallara:2000}). Sets of this form are called \textit{sets of finite perimeter}. As a consequence of Rellich's theorem in $BV$ and the lower semicontinuity of the perimeter we get the following result (see \cite[Proposition 3.38]{Ambrosio-Fusco-Pallara:2000})).
\begin{theorem}\label{clea-th: set per}
Let $\Omega \subset \R^d$ open, bounded. Let $(W_k)_k \subset \Omega$ be a sequence of measurable sets with ${\cal H}^{d-1}(\partial^* W_k) \le C$ for some constant $C$ independent of $k$. Then there is a subsequence (not relabeled) and a measurable set $W$ such that $\chi_{W_k} \to \chi_W$ in measure for $k \to \infty$ and ${\cal H}^{d-1}(\partial^* W) \le \liminf_{k \to \infty} {\cal H}^{d-1}(\partial^* W_k)$.
\end{theorem}

We say that a function $u \in L^1(\Omega, \R^d)$ is in $BD(\Omega)$ if the  symmetrized distributional derivative $Eu := \frac{1}{2}((Du)^T + Du)$ 
is a finite $\R^{d \times d}_{\rm sym}$-valued Radon measure. Likewise, we say $u$ is  a \emph{special  function of bounded deformation} if $Eu$ has vanishing Cantor part $E^c u$. Then $Eu$ can be decomposed as 
\begin{align}\label{rig-eq: symmeas}
 Eu = e(u) {\cal L}^d  + E^s u = e(u) {\cal L}^d + [u] \odot \xi_u {\cal H}^{d-1}|_{J_u},
 \end{align}
where $e(u)$ is the absolutely continuous part of $Eu$ with respect to the Lebesgue measure ${\cal L}^d$, $[u]$, $\xi_u$, $J_u$ as before and $a \odot b = \frac{1}{2}(a \otimes b + b \otimes a)$. If in addition $e(u) \in L^p(\Omega)$ for $1 < p <\infty$ and ${\cal H}^{d-1}(J_u) < \infty$, we write $u \in SBD^p(\Omega)$. For basic properties of this function space we refer to \cite{Ambrosio-Coscia-Dal Maso:1997,  Bellettini-Coscia-DalMaso:98}.

 By ${\cal R}(\Omega) = \lbrace T: \Omega \to \R^d: T(x) = A\,x+b, A\in \R^{d\times d}_{\rm skew}, b \in \R^d\rbrace$ we denote  the space of infinitesimal rigid motions.  We recall a Korn-Poincar\'e inequality in $BD$ (see \cite{Kohn:82, Temam:85}). 

\begin{theorem}\label{rig-theo: korn-poini}
Let $\Omega \subset \R^d$ open, bounded,  connected with Lipschitz boundary.  Then there exist  a  linear  continuous map $P: BD(\Omega) \to {\cal R}(\Omega)$ which leaves the elements of ${\cal R}(\Omega)$ fixed and   a constant $C=C(\Omega)>0$ such that for all $u \in BD(\Omega)$
$$\Vert u - Pu \Vert_{L^{\frac{d}{d-1}}(\Omega)} \le C |Eu|(\Omega),$$
where $|Eu|(\Omega)$ denotes the total variation of $Eu$. The constant $C$ is invariant under rescaling of the domain.
\end{theorem}

\subsection{Poincar\'e's and Korn's inequality}\label{rig-sec: korn}

A key idea in our analysis will be the replacement of displacement fields in $SBD$ by suitable Sobolev functions and then the application of well know Poincar\'e and Korn inequalities. As the estimates will be employed on different Lipschitz sets, we need to provide uniform bounds for the constants involved in the inequalities. To this end, we introduce the notion of \emph{John domains}.

\begin{definition}\label{def: chain} 
{\normalfont

Let $\Omega \subset \R^d$ be an open, bounded set and let $x_0 \in \Omega$. We say $\Omega$ is a $c$-\emph{John domain} with respect to the \emph{John center} $x_0$ and with the constant
$c$ if for all $x \in \Omega$ there exists a rectifiable curve $\gamma: [0,l_\gamma] \to \Omega$, parametrized by arclength, such that $\gamma(0) = x$, $\gamma(l_\gamma) = x_0$ and 
$t \le c\dist(\gamma(t),\partial \Omega)$ for all $t \in [0,l_\gamma]$. 
}
 \end{definition}
 
 Domains of this form were introduced by John in \cite{John:1961} to study problems in elasticity theory and the term was first used by Martio and Sarvas in  \cite{Martio-Sarvas:1978}. Roughly speaking, a domain is a John domain if it is possible to connect two arbitrary points without getting too close to the boundary of the set. This class is much larger than Lipschitz domains and contains sets which may possess fractal boundaries or internal cusps (external cusps are excluded), e.g. the interior of Koch's \emph{snow flake} is a John domain.

Although in the present work only Lipschitz sets occur, it is convenient to consider this more general notion as the constants in  Poincar\'e's and Korn's inequalities only depend on the John constant. More precisely, we have the following statement (see e.g. \cite{Acosta, BuckleyKoskela, Diening}).

\begin{theorem}\label{th: kornsobo}
Let ­ $\Omega \subset \R^d$ be a $c$-John domain. Let $p \in (1,\infty)$ and $q \in (1,d)$. Then there is a constant $C=C(c,p,q)>0$ such that for all $u \in W^{1,p}(\Omega)$ there is some $A \in \R^{d \times d}_{\rm skew}$ such that
$$\Vert \nabla u -A\Vert_{L^p(\Omega)} \le C\Vert e(u) \Vert_{L^p(\Omega)}.$$
Moreover,  for all $u \in W^{1,q}(\Omega)$ there is some $b \in \R^d$ such that
$$\Vert u -b\Vert_{L^{q^*}(\Omega)} \le C\Vert \nabla u \Vert_{L^q(\Omega)},$$
where $q^* = \frac{dq}{d-q}$. The constant is invariant under rescaling of the domain. 
\end{theorem}

\subsection{A modification with controllable jump heights}\label{rig-sec: modifica}

A main strategy of our proof will be the application of Theorem \ref{rig-theo: korn-poini} in certain regions of the domain. It first appears that this inequality is not adapted for the estimates in \eqref{eq: main estmain}  as in $|Eu|(\Omega)$ not only the elastic but also the surface energy depending on the jump height is involved. However, in \cite{Friedrich:15-1} we have shown that one can indeed find bounds on the jump heights in terms of the elastic energy  after a suitable modification of the displacement field.

In the following $d(W)$ denotes the diameter of a set $W \subset \R^2$ and for $\mu>0$ we define $Q_\mu:=(-\mu,\mu)^2$. We recall the following result (see \cite[Theorem 2.2]{Friedrich:15-1}).

\begin{theorem}\label{rig-th: derive prop}
Let $\lambda> 0$. Then  there is a constant $C=C(\lambda)$  and a universal constant $c>0$ both independent of $\mu$ such that for all $\eps>0$, $\delta >0$ and all  $u \in SBD^2(Q_\mu)$ the following holds: There are paraxial rectangles $R_1,\ldots,R_n$ with 
\begin{align}\label{eq:ers3}
\sum\nolimits^n_{j=1} d(R_j) \le (1 + c \lambda)\big( {\cal H}^1(J_u) + \eps^{-1}\Vert e(u) \Vert^2_{L^2(Q_\mu)}\big)
\end{align}
and a modification $\bar{u} \in SBV^2( Q_{\tilde{\mu}}  )$ with $J_{\bar{u}} \subset \bigcup^n_{j=1} \partial R_j$ and 
\begin{align}\label{eq:erst5}
\Vert \bar{u}- u\Vert_{L^2(Q_{\tilde{\mu}} \setminus G)} \le \delta, \ \  \ \ \Vert e(\bar{u})\Vert_{L^2(Q_{\tilde{\mu}})} \le  \Vert e(u)\Vert_{L^2(Q_\mu)} + \delta,
\end{align}
 where $G :=\bigcup^n_{j=1} R_j$ and $\tilde{\mu} = \max\lbrace \mu - 3\sum_j d(R_j), 0\rbrace$, such that for all measurable sets $D \subset Q_{\tilde{\mu}}$ we have
\begin{align}\label{eq:erst4}
(|E \bar{u}|(D))^2  \le c |D| \Vert e(\bar{u}) \Vert^2_{L^2( D )} +  C\eps {\cal H}^1(D \cap J_{\bar{u}} ) \sum\nolimits_{R_j \in {\cal R}(D)}(d(R_j))^2,
\end{align}
where ${\cal R}(D):= \lbrace R_j: D \cap \overline{R_j} \neq \emptyset \rbrace$. Moreover, if $u \in L^\infty(Q_\mu)$, we can choose the modification such that $\Vert \bar{u} \Vert_{L^\infty(Q_{\tilde{\mu}} \setminus G)} \le \Vert u \Vert_{L^\infty(Q_\mu)}$.
\end{theorem}

\section{The local estimate on a square}\label{sec: local}

This section is devoted to the derivation of a local estimate on a square. Recall $Q_\mu = (-\mu,\mu)^2$ for $\mu >0$. 

\begin{theorem}\label{th: korn1}
Let $p \in [1,2), q \in [1,\infty)$ and $\mu>0$. Then there is a constant $C=C(p,q)>0$ independent of $\mu$ such that for all $u \in SBD^2(Q_\mu)$ there is a set of finite perimeter $F \subset Q_\mu$ with 
\begin{align}\label{eq: R2}
{\cal H}^1(\partial^* F) \le C{\cal H}^1(J_u), \ \ \ \ |F| \le C({\cal H}^1(J_u))^2
\end{align}
and $A \in \R^{2 \times 2}_{\rm skew}$, $b \in \R^2$ such that
\begin{align}\label{eq: main est}
\begin{split} 
(i) & \ \ \Vert u - (A\,\cdot + b) \Vert_{L^q(Q_{\bar{\mu}} \setminus F)}\le C \mu^{\frac{2}{q}} \Vert e(u) \Vert_{L^2( Q_\mu)},\\
(ii) & \ \ \Vert \nabla u - A \Vert_{L^p(Q_{\bar{\mu}} \setminus F)}\le C\mu^{\frac{2}{p}-1} \Vert e(u) \Vert_{L^2( Q_\mu)},
\end{split}
\end{align}
where $Q_{\bar{\mu}} = (-\bar{\mu},\bar{\mu})^2$ with $\bar{\mu} = \max\lbrace \mu - C{\cal H}^1(J_u) ,0 \rbrace$. Moreover, we obtain 
$u \chi_{Q_{\bar{\mu}} \setminus F} \in GSBV^p(Q_\mu)$. If in addition $u \in L^\infty(Q_\mu)$, we find $u \chi_{Q_{\bar{\mu}} \setminus F} \in SBV^p(Q_\mu)$.  
\end{theorem}  

Observe that the additional statement that $u \chi_{Q_{\bar{\mu}} \setminus F} $ lies in $(G)SBV^p$ does not directly follow from \eqref{eq: main est}(ii) as a Korn-type inequality 
for $\nabla u$ does not automatically guarantee that $u$ has bounded variation. In fact, the property $|D^c u|(Q_{\bar{\mu}}) \le \sqrt{2}|E^cu|(Q_{\bar{\mu}})$, which holds for $BV$ functions (see \cite{Alberti:93}), is not known in $BD$.  (We remark that the analog of Alberti's rank one property in $BD$   \cite{rindler} is not expedient here.)   
 In the present context, we circumvent this problem by  an approximation of $SBD$ functions.

We first show the result for modifications given by Theorem \ref{rig-th: derive prop} and afterwards we  prove the general version of Theorem \ref{th: korn1} by considering sequences of modifications.

\subsection{Local estimate for modifications}

We first introduce some further notation. For $s>0$ we partition $\R^2$ up to a set of measure zero into squares $Q^s(p) = p + s(-1,1)^2$ for $p \in I^s := s(1, 1) + 2s\Z^2$.   Let $\theta \in 2^{-\N}$ small, fixed and define $s_i = \mu \theta^i$ for $i \ge 0$. We let
\begin{align}\label{eq: enlarged squares}
{\cal Q}_i = \lbrace Q=Q^{s_i}(p): p \in I^{s_i} \rbrace.
\end{align}
For each $Q \in {\cal Q}_i$ we introduce enlarged squares $Q  \subset Q'' \subset Q'$ defined by
\begin{align}\label{eq:enlarged}
Q'' = \tfrac{5}{4}Q, \ \ \ Q' = \tfrac{3}{2}Q,
\end{align}
where $\lambda Q$ denotes the square with the same center and $\lambda$-times the sidelength of $Q$.
Moreover, by $\dist(A,B)$ we denote the  Euclidian  distance between $A,B \subset \R^2$.  In the sequel, \emph{infinitesimal rigid motions} $A\,x + b$, $A\in \R^{2 \times 2}_{\rm skew}$, $b \in \R^2$, will appear frequently and we will often write $a(x)=a_{A,b}(x) = A\,x + b$ for the sake of brevity. We now prove the local result for modifications.

\smallskip

\begin{theorem}\label{lemma: korn2}
Let $p \in [1,2), q \in [1,\infty)$ and $\mu>0$.  Then there is a constant $C=C(p,q)>0$ such that for each $\delta >0$ and  all $u \in  SBD^2(Q_\mu)$  there is a finite union of rectangles $G \subset Q_\mu$ with ${\cal H}^1(\partial G) \le C{\cal H}^1(J_u)$ and a modification $\bar{u} \in SBV^2(Q_{\bar{\mu}})$ for  $\bar{\mu} = \max\lbrace \mu - C{\cal H}^1(J_u) ,0 \rbrace$ with
\begin{align}\label{eq:modi}
\begin{split}
(i)& \ \ \Vert \bar{u}- u\Vert_{L^2(Q_{\bar{\mu}} \setminus G)} \le \delta, \ \  \ \    \Vert e(\bar{u})\Vert_{L^2(Q_{\bar{\mu}})} \le  \Vert e(u)\Vert_{L^2(Q_\mu)} + \delta, \\ 
(ii)& \ \  {\cal H}^1(J_{\bar{u}}) \le C{\cal H}^1(J_u)
\end{split}
\end{align}
such that there is a set $F \subset Q_\mu$ satisfying \eqref{eq: R2} and $A \in \R^{2\times 2}_{\rm skew}$, $b \in \R^2$ with
\begin{align}\label{eq:modi2}
\mu^{-\frac{2}{q}}\Vert  \bar{u} - (A\,\cdot + b) \Vert_{L^q(Q_{\bar{\mu}} \setminus F)} + \mu^{1-\frac{2}{p}}\Vert \nabla \bar{u} - A \Vert_{L^p(Q_{\bar{\mu}} \setminus F)}\le C (\Vert e(u) \Vert_{L^2( Q_\mu)} + \delta).
\end{align}  
\end{theorem}

\Proof Let $u \in  SBD^2(Q_\mu)$ and $\delta>0$ be given. Without restriction we can assume $\mathcal{H}^1(J_u)>0$ as otherwise the statement follows from Theorem \ref{th: kornsobo}. We first apply Theorem \ref{rig-th: derive prop} with $\eps = ({\cal H}^1(J_u))^{-1}    \Vert e(u) \Vert^2_{L^2(Q_\mu)}$ and $\lambda = 1$ to obtain a modification $\bar{u} \in SBV^2(\tilde{Q})$ such that \eqref{eq:ers3}-\eqref{eq:erst4} hold with  $ \tilde{Q} := Q_{\tilde{\mu}}  \subset Q_\mu$ as in Theorem \ref{rig-th: derive prop}. Let $\bar{J}:= \bigcup_{l=1}^n \partial R_l$ and $G:= \bigcup^n_{l=1} R_l$. We particularly have by \eqref{eq:ers3}
\begin{align}\label{eq: R1}
{\cal H}^1(\partial G) \le {\cal H}^1(\bar{J}) \le 2\sqrt{2}\sum\nolimits^n_{l=1} \diam(R_l) \le 4\sqrt{2}(1+c){\cal H}^1(J_u).
\end{align}
Choosing the constant $C>0$ in the assertion large enough (depending on $\theta$) we get  by \eqref{eq: R1}  
\begin{align}\label{eq: new1}
Q_{\bar{\mu}}\subset \tilde{Q}, \ \ \ \  \  \ \  \sum\nolimits_{l=1}^n \diam(R_l) \le \tfrac{\theta}{24}\dist(\partial \tilde{Q}, \partial Q_{\bar{\mu}}).
\end{align}
Then in view of \eqref{eq:erst5}, \eqref{eq: R1} and $J_{\bar{u}} \subset \bar{J}$, we get  \eqref{eq:modi}. 

We now derive \eqref{eq:modi2} for the modification $\bar{u}$, where we will regard $\bar{u}$ as a function defined on $\tilde{Q}$.  We may assume that $Q_{\bar{\mu}} \neq \emptyset$ as otherwise the assertion of the lemma is trivial. Then \eqref{eq: new1} implies 
\begin{align}\label{NNN}
\sum\nolimits_l \diam(R_l) \le \tfrac{\theta}{24} \mu = \tfrac{1}{24} s_1.
\end{align}
We start with the identification of regions where  the ${\cal H}^1$ measure of  $\bar{J}$ is too large (Step I). Afterwards, we will use \eqref{eq:erst4}  to apply Theorem \ref{rig-theo: korn-poini} on these specific sets (Step II). This Korn-Poincar\'e estimate will then enable us to define a suitable modification (Step III) for which Korn's inequality for Sobolev functions can be used (Step IV). We also refer to Figure \ref{bild2} below, where some of the objects introduced in the proof are illustrated.

\smallskip

\noindent  {\bf Step I (Identification of `bad' sets):} We first identify squares of various length scales where  the ${\cal H}^1$ measure of  $\bar{J}$ is too large. Recalling  \eqref{eq: enlarged squares} and \eqref{eq:enlarged} we introduce the sets 
\begin{align}\label{eq: calA} 
{\cal A}_i = \big\{ Q \in {\cal Q}_i: \ Q'' \subset \tilde{Q}, \ \   \sum\nolimits_{l=1}^n{\cal H}^1(Q'' \cap \partial R_l) \ge \tfrac{1}{8}\theta  s_i \big\}
\end{align}
for $i \in \N$. Let $A_i = \bigcup_{Q \in {\cal A}_i} \overline{Q}$. Then we define $${\cal B}_i = \big\lbrace Q \in {\cal A}_i: \ Q\cap Q_{\bar{\mu}} \neq \emptyset, \  \ Q \cap \bigcup\nolimits^{i-1}_{j=1} A_j = \emptyset \big\rbrace$$
and accordingly let  $B_i = \bigcup_{Q \in {\cal B}_i} \overline{Q}$ and $B'_i = \bigcup_{Q \in {\cal B}_i} \overline{Q'}$ for $i \in \N$. For later we note that $\bigcup_{i \ge 1} A_i \cap Q_{\bar{\mu}}= \bigcup_{i \ge 1} B_i \cap Q_{\bar{\mu}}$. We now show that for some $I \in \N$ sufficiently large we have
\begin{align}\label{eq: prop1}
\begin{split}
(i) & \ \ \sum\nolimits_{l: \, R_l \cap Q'' \neq \emptyset} \diam(R_l) \le  s_i  \ \ \text{ for all } \ Q \in {\cal B}_i, \, i\ge 1,\\
(ii) & \ \ B_i = \emptyset \ \text{ for } \ i> I,\\
(iii) & \ \ J_{\bar{u}}  \cap Q_{\bar{\mu}} \subset \bar{J}  \cap Q_{\bar{\mu}}\subset \bigcup\nolimits_{i \ge 1} B_i,\\
(iv) & \ \ {\cal H}^1\big(\partial \, \big(\bigcup\nolimits_{i\ge 1}B'_i \big)\big) \le C\theta^{-1}\sum\nolimits_{l=1}^n d(R_l) \le C\mu,\\
(v) & \ \ |B'_i| \le Cs_i \theta^{-1}\sum\nolimits_{l=1}^n d(R_l)\le Cs_i\mu \ \text{ for all } i\ge 1
\end{split}
\end{align}
for a constant $C>0$ independent of $\theta$. 

We first confirm (i).  We assume $\sum\nolimits_{R_l \in {\cal F}} \diam(R_l) >  s_i$ for some $Q \in {\cal B}_i$, where ${\cal F} = \lbrace R_l: R_l \cap Q'' \neq \emptyset\rbrace$,  and derive a contradiction. Choose $j \le i$ such that  $\frac{1}{8}s_j < \sum\nolimits_{R_l \in {\cal F}} \diam(R_l) \le \frac{1}{8}s_{j-1}$ and observe that by  \eqref{NNN}  we find $j \ge 2$.   Moreover, we select $Q_* \in {\cal Q}_{j-1}$ such that $Q \subset Q_*$. As $s_{j-1} < 8\theta^{-1} \sum_l \diam(R_l)  \le \frac{1}{3}\dist(\partial \tilde{Q}, \partial Q_{\bar{\mu}})$   by \eqref{eq: new1} and $Q_* \cap Q_{\bar{\mu}} \neq \emptyset$, we find $Q''_* \subset \tilde{Q}$ recalling \eqref{eq:enlarged}.  Since $R_l \cap Q'' \neq \emptyset$,   $d(R_l) \le \frac{1}{8}s_{j-1}$ for all $R_l \in {\cal F}$ and  $Q \subset Q_*$,   we get   $\overline{R_l} \subset Q_*''$ for all $R_l \in {\cal F}$.  Thus,  
$$\sum\nolimits_{l=1}^n{\cal H}^1(\partial R_l \cap Q''_*) \ge \sum\nolimits_{R_l \in {\cal F}}{\cal H}^1(\partial R_l)    > \tfrac{1}{8} s_{j} = \tfrac{1}{8}\theta s_{j-1}.$$ 
This yields $Q_* \in {\cal A}_{j-1}$. Consequently, as $Q \subset Q_*$, this implies $Q \notin {\cal B}_i$ giving the desired contradiction.

 To see  (ii), choose $I \in \N$ so large that $s_{I+1} < \min_{l=1}^n d(R_l)$. Assume there was some $Q \in {\cal B}_i$ for $i >I$. Then we find some $R_k$ with $R_k \cap Q'' \neq \emptyset$ and $\diam(R_k)  > s_{I+1} \ge s_{i}$.  This implies $\sum\nolimits_{l: \, R_l \cap Q'' \neq \emptyset} \diam(R_l) > s_{i}$ and yields a contradiction to \eqref{eq: prop1}(i).  

Moreover, the definition in  \eqref{eq: calA} implies $\bar{J}  \cap Q_{\bar{\mu}} \subset \bigcup_{i \ge 1} A_i \cap Q_{\bar{\mu}}$  and thus (iii) follows from the property $\bigcup_{i \ge 1} A_i \cap Q_{\bar{\mu}}= \bigcup_{i \ge 1} B_i \cap Q_{\bar{\mu}}$ and the fact that $J_{\bar{u}} \subset \bar{J}$. 

To show (iv), we define $\hat{\cal B}_i = \lbrace Q \in {\cal B}_i: \overline{Q'} \not\subset \bigcup^{i-1}_{j=1} B'_j\rbrace$ as well as $\hat{B}'_i = \bigcup_{Q \in \hat{\cal B}_i} \overline{Q'}$ and observe $\bigcup_{i\ge 1} B_i' = \bigcup_{i\ge 1} \hat{B}_i'$. It is elementary to see that for $\theta$ small enough each $x \in \tilde{Q}$ is contained in at most one  set $\bigcup_{Q \in \hat{\cal B}_i} \overline{Q''}, i \in \N$, and thus contained in at most four different squares $Q''$ with $Q \in \bigcup_{i \ge 1} \hat{\cal B}_i$. Consequently, we derive using  \eqref{eq: calA} for a universal constant $C>0$ large enough
\begin{align*}
{\cal H}^1\big(\partial \, \big(\bigcup\nolimits_{i\ge 1}B'_i \big)\big)  &\le \sum_{Q \in \bigcup_{i\ge 1}\hat{\cal B}_i}{\cal H}^1(\partial Q')\le C\theta^{-1}\sum_{Q \in \bigcup_{i\ge 1}\hat{\cal B}_i} \Big(\sum\nolimits_{l=1}^n{\cal H}^1(Q'' \cap \partial R_l) \Big) \\
&\le C\theta^{-1}\sum\nolimits_{l=1}^n{\cal H}^1(\tilde{Q} \cap \partial R_l) \le C\theta^{-1}\sum\nolimits_{l=1}^n d(R_l) \le C\mu,
\end{align*}
where in the last step we employed \eqref{eq: new1}. Finally, (v) follows from a similar argumentation, again using \eqref{eq: calA}.

\vspace{-0.2cm}
\begin{figure}[H]
\centering
\begin{overpic}[width=1.0\linewidth,clip]{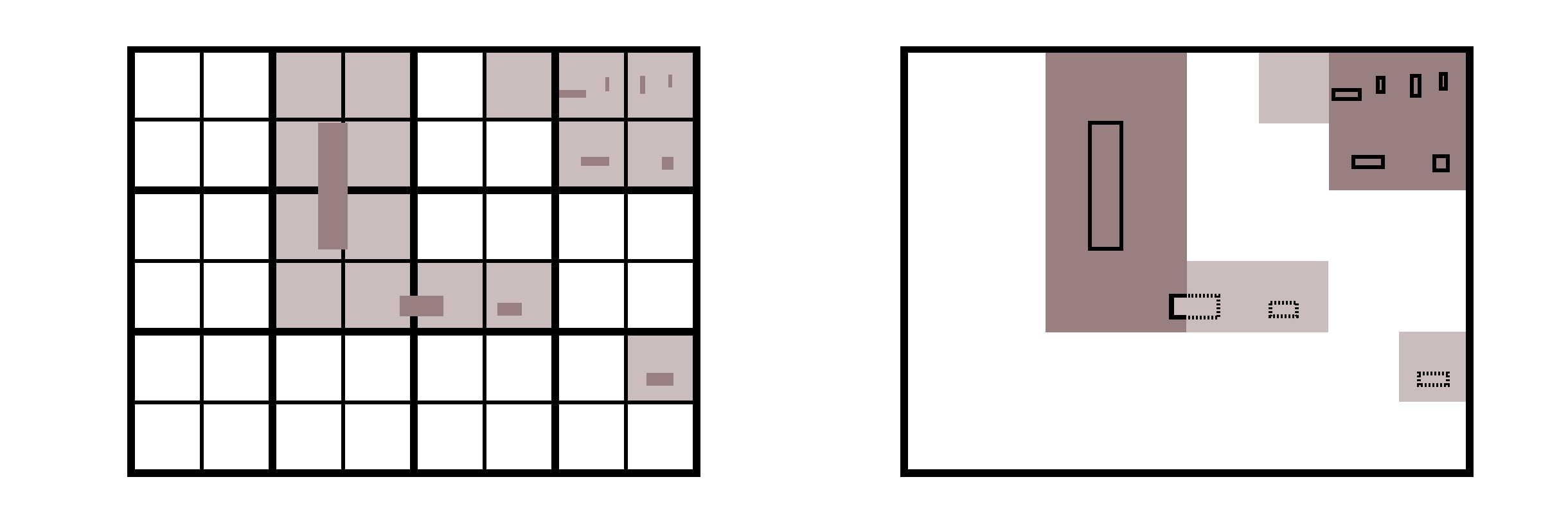}
\put(16,116){{$(a)$}}
\put(36,12){{$\underbrace{\ \ \ \ \ \ \ \ \   }$}}
\put(60,0){{\small $s_{i-1}$}}
\put(90,12){{\footnotesize $\underbrace{\  }$}}
\put(100,0){{\small $s_{i}$}}
\put(220,116){{$(b)$}}
\put(246,56){{$B_{i-1}$}}
\put(348,36){{$B_{i}$}}
\put(260,86){$\bar{J}$}
\put(270,60){\line(1,0){12}}
\put(360,43){\line(1,0){12}}
\put(270,90){\line(1,0){17}}
\end{overpic}
\caption{\small  (a) Sketch of a part of $Q_\mu$ with the squares ${\cal Q}_{i-1}$ and ${\cal Q}_i$ for $\theta =\frac{1}{2}$. In dark gray the union of rectangles $G$ is depicted and in light gray the squares contained in $A_{i-1} \cup A_i$. (b) The sets $B_{i-1}$ (dark gray) and $B_i$ (light gray) are illustrated, where $\bar{J} \subset B_{i-1} \cup B_i$. Note that the ${\cal H}^1$ measure of $\bar{J}$ in the (enlarged) squares of ${\cal B}_{i-1}, {\cal B}_i$ is not `too large' such that \eqref{eq: prop1}(i) holds and also not `too small' (see \eqref{eq: calA}) such that $|B'_{i-1}|, |B'_i|$ is suitably controlled, cf. \eqref{eq: prop1}(v). In the modifcation $\bar{u}_{I-i+1}$ the discontinuities of $\bar{u}$ depicted in dashed lines are removed, in $\bar{u}_{I-i+2}$ the remaining part of $\bar{J}$ is removed.  } \label{bild2}
\end{figure}
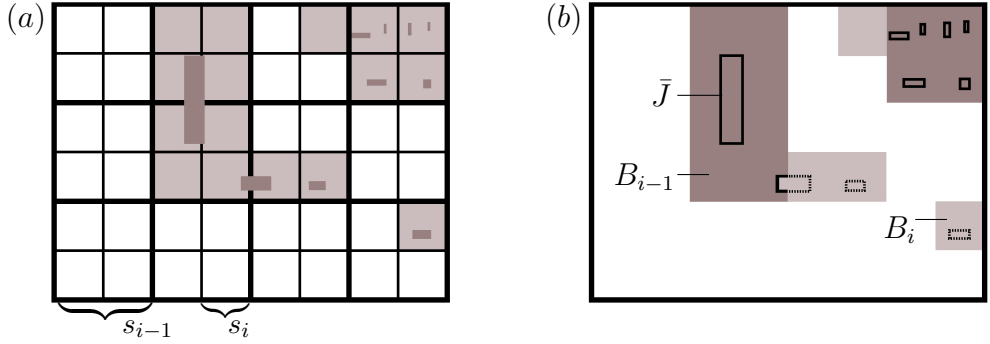

\noindent {\bf Step II (Korn-Poincar\'e inequality):} Recall that by the definition in \eqref{eq: calA}  we find $Q'' \subset \tilde{Q}$ for every $Q \in {\cal B}_i$. Thus,  using Theorem \ref{rig-theo: korn-poini} on $Q''$ for $Q \in {\cal B}_i$ we obtain by \eqref{eq:erst4} and \eqref{eq: prop1}(i) infinitesimal rigid motions $a_Q = a_{A_Q,b_Q}$ such that 
\begin{align*}
\Vert \bar{u} - a_Q\Vert^2_{L^2(Q'')} &\le C (|E\bar{u}|(Q''))^2  \\&\le Cs^2_i \Vert e(\bar{u}) \Vert^2_{L^2(Q'')} +  C\eps{\cal H}^1(J_{\bar{u}} \cap Q'')\sum\nolimits_{R_l \cap Q'' \neq \emptyset} (\diam(R_l))^2 \\
&\le Cs^2_i \big(\Vert e(\bar {u}) \Vert^2_{L^2(Q'')} + \eps{\cal H}^1(\bar{J} \cap Q'') \big)
\end{align*}
for all $Q \in {\cal B}_i$, where in the last step we employed $J_{\bar{u}} \subset \bar{J}$. For shorthand, we let $ {\cal E}  : = \Vert e(u) \Vert_{L^2(Q_\mu)} + \delta$. Summing over all squares in ${\cal B}_i$, using \eqref{eq:erst5} and recalling that each $x\in \tilde{Q}$ is contained in at most four  different $Q''$, $Q \in {\cal B}_i$, we get for all $i \ge 1$
\begin{align}\label{eq: sum}
\sum\nolimits_{Q \in {\cal B}_i} \Vert \bar{u} - a_Q\Vert^2_{L^2(Q'')} \le Cs_i^2  (\Vert e(\bar{u}) \Vert^2_{L^2(Q_\mu)} + \Vert e(u) \Vert^2_{L^2(Q_\mu)}) \le Cs_i^2{\cal E}^2,
\end{align}
where we used  $\eps = ({\cal H}^1(J_u))^{-1}    \Vert e(u) \Vert^2_{L^2(Q_\mu)}\le C ({\cal H}^1(\bar{J}))^{-1} \Vert e(u) \Vert^2_{L^2(Q_\mu)}$ (cf. \eqref{eq: R1}).

\smallskip 

\noindent {\bf Step III (Modification):} We now show that we can `heal' the discontinuities of $\bar{u}$ in $Q_{\bar{\mu}}$. The strategy is to modify the displacement field inductively. Let $1 < p <2$ be given. Let $I \in \N$ be the largest index such that $B_{I} \neq \emptyset$ (see \eqref{eq: prop1}(ii)). Define $\bar{u}_0 = \bar{u}$ and assume $\bar{u}_j \in SBD^2(\tilde{Q})$ has already been constructed satisfying
\begin{align}\label{eq: prop2}
(i) & \  \ J_{\bar{u}_j} \cap Q_{\bar{\mu}} \subset \bigcup\nolimits^{I-j}_{k=1} B_k,\\
(ii) & \ \ \sum\nolimits_{Q \in {\cal B}_i} \Vert \bar{u}_j - a_Q\Vert^2_{L^2(Q'')} \le \bar{C} s_i^2  \,  \prod\nolimits^j_{k=0} (1+ \eta^{I-i-k}) \  {\cal E}^2  \ \,  \forall i \le I-j,  \notag\\
(iii) & \ \ \Vert e(\bar{u}_j)\Vert^p_{L^p(Q_{\bar{\mu}})}  \le  \bar{C}\mu^{2-p} \, \prod\nolimits^j_{k=0} (1+\eta^{I-k}) \  {\cal E}^p \notag
\end{align}
for some $\bar{C}$ large enough, where for shorthand $\eta=\theta^{\frac{1}{2}-\frac{p}{4}}<1$. Clearly, by \eqref{eq:modi}, \eqref{eq: prop1}(ii),(iii), \eqref{eq: sum} and H\"older's inequality together with $|Q_{\bar{\mu}}|\le 4\mu^2$ we find that \eqref{eq: prop2} holds for $j=0$.

We now construct $\bar{u}_{j+1}$. In the following $C>0$ denotes a generic constant which is always independent of $\theta$. We consider a partition of unity $\lbrace \varphi_0 \rbrace \cup (\varphi_Q)_{Q \in {\cal B}_{I-j}} \subset C^\infty(\R^2)$ with the properties 
\begin{align}\label{eq: prop3}
\begin{split}
(i) & \ \ \varphi_0(x) + \sum\nolimits_{Q \in {\cal B}_{I-j}} \varphi_Q(x) = 1 \ \text{ for all }x \in  \tilde{Q},\\ 
(ii) & \  \ Q \subset {\rm supp}(\varphi_Q) \subset Q'' \ \text{ for all } Q \in {\cal B}_{I-j},\\
(iii) & \ \ {\rm supp}(\varphi_0) \subset \R^2 \setminus B_{I-j},\\
(iv) & \ \ \Vert \nabla \varphi_0\Vert_\infty, \Vert \nabla \varphi_Q \Vert_\infty \le Cs_{I-j}^{-1} \ \text{ for all } \ Q \in {\cal B}_{I-j}.
\end{split}
\end{align}
Then we define
$$ \bar{u}_{j+1}(x) = \bar{u}_j(x) + \sum_{Q \in {\cal B}_{I-j}} \varphi_Q(x) (A_Q \, x + b_Q - \bar{u}_j(x)) = \bar{u}_j + \sum_{Q \in {\cal B}_{I-j}} \varphi_Q (a_Q - \bar{u}_j)$$
for all $x \in \tilde{Q}$.  As $\lbrace x\in \tilde{Q}: \sum\nolimits_{Q \in {\cal B}_{I-j}} \varphi_Q(x) = 1 \rbrace \supset B_{I-j} $  by \eqref{eq: prop3}(i),(iii), we get $\bar{u}_{j+1} = \sum_{Q \in {\cal B}_{I-j}} \varphi_Q a_Q$ in $B_{I-j}$. Thus,   $\bar{u}_{j+1}$ is smooth in $B_{I-j}$  and \eqref{eq: prop2}(i) holds. Using \eqref{eq: prop2}(ii) for $j$ and $i = I-j$ we obtain 
\begin{align*}
 \Vert \bar{u}_{j+1} - \bar{u}_j\Vert^2_{L^2(\tilde{Q})} &\le C\sum\nolimits_{Q \in {\cal B}_{I-j}}\Vert \bar{u}_j - a_Q\Vert^2_{L^2(Q'')} \\
& \le C\bar{C} s^2_{I-j} \,  \prod\nolimits^j_{k=0} (1+ \eta^{j-k}) {\cal E}^2 \le C\bar{C} s^2_{I-j}{\cal E}^2,
\end{align*}
where in the first step we used that each $x \in \tilde{Q}$ is contained in at most four different enlarged squares. Using $\eta  = \theta^{\frac{1}{2}-\frac{p}{4}} \ge \theta$, we get for $\theta$ small enough (recall that $p \in (1,2)$ and that  $C$ is independent of $\theta$) 
$$\sum\nolimits_{Q \in {\cal B}_i}\Vert \bar{u}_{j+1} - \bar{u}_j\Vert^2_{L^2(Q'')} \le C\bar{C} s^2_{I-j}{\cal E}^2  \le C\bar{C} s_i^2\eta^{2(I-j-i)}{\cal E}^2 \le \tfrac{1}{6}\bar{C}s_i^2 \eta^{2(I-j-i-1)}{\cal E}^2$$
for all $1 \le i \le I-j-1$. The previous estimate together with \eqref{eq: prop2}(ii) for $j$  and a scaled version of Young's inequality of the form $(a+b)^2 \le  (1+ \delta)a^2 + (1+\frac{1}{\delta}) b^2$ ($a,b  \in \R$, $\delta >0$) yields   for $\delta = \tfrac{1}{2} \eta^{I-i-j-1}$
\begin{align*}
&\sum\nolimits_{Q \in {\cal B}_i} \Vert \bar{u}_{j+1} - a_Q\Vert^2_{L^2(Q'')} \\
& \ \ \ \ \ \le \sum\nolimits_{Q \in {\cal B}_i} \Big((1+\tfrac{1}{2} \eta^{I-i-j-1}) \Vert \bar{u}_{j} - a_Q\Vert^2_{L^2(Q'')} +  3  \eta^{-(I-i-j-1)}\Vert \bar{u}_{j+1} - \bar{u}_j\Vert^2_{L^2(Q'' )}\Big)\\
& \ \ \ \ \ \le \bar{C} s_i^2  \,  \prod\nolimits^j_{k=0} (1+ \eta^{I-i-k})(1+\tfrac{1}{2} \eta^{I-i-j-1}){\cal E}^2  + \bar{C} \tfrac{1}{2} \eta^{I-i-j-1} s^2_{i} {\cal E}^2 \\ &\ \ \ \ \
 \le \bar{C} s_i^2  \,  \prod\nolimits^{j+1}_{k=0} (1+ \eta^{I-i-k}) \, {\cal E}^2
\end{align*}
for all $1 \le i \le I-j-1$. This shows \eqref{eq: prop2}(ii). To confirm  \eqref{eq: prop2}(iii), we first note that by H\"older's inequality, \eqref{eq: prop2}(ii) and the fact that $|B'_{I-j}| \le C \mu s_{I-j}$ (see \eqref{eq: prop1}(v)) we obtain
\begin{align}\label{eq: once}
\begin{split}
\sum\nolimits_{Q\in{\cal B}_{I-j}} \Vert \bar{u}_{j} - a_Q\Vert^p_{L^p(Q'')} &\le C|B'_{I-j}|^{1-\frac{p}{2}}\big(\sum\nolimits_{Q\in{\cal B}_{I-j}} \Vert \bar{u}_{j} - a_Q\Vert^2_{L^2(Q'')}\big)^\frac{p}{2}\\
& \le C\bar{C}^{\frac{p}{2}}\mu^{1-\frac{p}{2}} s^{1+\frac{p}{2}}_{I-j} {\cal E}^p,
\end{split}
\end{align}
where we again used that each $x\in \tilde{Q}$ is contained in at most four different enlarged squares. We calculate the derivative
$$\nabla \bar{u}_{j+1}  =  \nabla \bar{u}_{j} \varphi_0 + \sum\nolimits_{Q \in {\cal  B}_{I-j}} \varphi_Q A_Q +   (a_Q - \bar{u}_j) \otimes \nabla \varphi_Q.$$
Now we  again apply a scaled version of Young's inequality of the form $|a+b|^p \le \big((1+ \delta)a^2 + (1+\frac{1}{\delta}) b^2\big)^{\frac{p}{2}} \le (1+ \delta^{\frac{p}{2}}) |a|^p + (1+\delta^{-\frac{p}{2}})|b|^p$ for $a,b  \in \R$, $\delta >0$. (Recall that $p \in (1,2)$.)  Consequently, similarly as before, using \eqref{eq: prop2}(iii),  \eqref{eq: prop3}(iv) and \eqref{eq: once}, we find with $\delta^{\frac{p}{2}} = \tfrac{1}{2}\eta^{I-j-1}$
\begin{align*}
\Vert e(\bar{u}_{j+1})\Vert^p_{L^p(Q_{\bar{\mu}})} &\le (1+\delta^{\frac{p}{2}})\Vert e(\bar{u}_{j})\Vert^p_{L^p(Q_{\bar{\mu}})} + C\delta^{-\frac{p}{2}}s_{I-j}^{-p} \sum\nolimits_{Q\in{\cal B}_{I-j}} \Vert \bar{u}_{j} - a_Q\Vert^p_{L^p(Q'')} \\ &\le
(1+\delta^{\frac{p}{2}})\Vert e(\bar{u}_{j})\Vert^p_{L^p(Q_{\bar{\mu}})} + C\delta^{-\frac{p}{2}}\bar{C}^{\frac{p}{2}}  s^{1-\frac{p}{2}}_{I-j} \mu^{1-\frac{p}{2}}{\cal E}^p\\
& \le \bar{C}\mu^{2-p} \, \prod\nolimits^j_{k=0} (1+\eta^{I-k}) \, (1+\tfrac{1}{2}\eta^{I-j-1}){\cal E}^p  + \bar{C} \tfrac{1}{2} \eta^{I-j-1} \mu^{2-p}{\cal E}^p \\&\le \bar{C}\mu^{2-p} \, \prod\nolimits^{j+1}_{k=0} (1+\eta^{I-k}){\cal E}^p,
\end{align*}
where we used $\eta^2=\theta^{1-\frac{p}{2}}$ and thus $Cs_{I-j}^{1-\frac{p}{2}} = C\mu^{1-\frac{p}{2}} \eta^{2(I-j)}\le \frac{1}{4} \mu^{1-\frac{p}{2}} \eta^{2(I-j-1)}$ for $\theta$ sufficiently small since $C$ is independent of $\theta$.

\smallskip 

\noindent {\bf Step IV (Korn's inequality):} Assume $\theta>0$ has been fixed according to Step III. We define $\hat{u} = \bar{u}_I$ and observe that by \eqref{eq: prop2}(i),(iii) we have $\hat{u}|_{Q_{\bar{\mu}}} \in W^{1,p}(Q_{\bar{\mu}})$ with 
\begin{align}\label{eq:modi3}
\Vert e(\hat{u}) \Vert_{L^p(Q_{\bar{\mu}})} \le C\mu^{\frac{2}{p}-1}{\cal E} = C\mu^{\frac{2}{p}-1}(\Vert e(u) \Vert_{L^2(Q_\mu)}+\delta)
\end{align}
for some $C=C(p)$. Moreover, we define $F = \bigcup^I_{i=1} B_i'$ and get that $\hat{u}  =\bar{u}$ on $Q_{\bar{\mu}} \setminus F$ due to the construction of the functions $(\bar{u}_j)_j$. By \eqref{eq: R1} and \eqref{eq: prop1}(iv)  we obtain ${\cal H}^1(\partial F) \le C{\cal H}^1(J_{{u}})$. 
In view of the definition in \eqref{eq: calA}, we find some $i_0 \in \N$ with  $s_{i_0} \le c\sum_{l=1}^n{\cal H}^1(\partial R_l \cap \tilde{Q})$ for a sufficiently large $c$ such that $B'_i = \emptyset$ for all $i\le i_0$. Thus, using \eqref{eq: R1}, \eqref{eq: prop1}(v) we find
$$\big|\bigcup\nolimits^I_{i=1} B_i'\big| \le C\sum\nolimits_l d(R_l)\,\sum\nolimits^I_{k=i_0+1} s_k \le C{\cal H}^1(J_u) s_{i_0} \le C({\cal H}^1(J_{u}))^2.$$   
This yields $|F| \le C({\cal H}^1(J_{u}))^2$ and shows \eqref{eq: R2}. We now apply Poincar\'e's and Korn's inequality (see Theorem \ref{th: kornsobo}) and find $A \in \R^{2 \times 2}_{\rm skew}$, $b \in \R^2$ such that by a standard scaling argument
\begin{align*}
 \Vert \nabla \bar{u} - A \Vert_{L^p(Q_{\bar{\mu}} \setminus F)} &\le \Vert \nabla \hat{u} - A \Vert_{L^p(Q_{\bar{\mu}})}\le C\Vert e(\hat{u}) \Vert_{L^p(Q_{\bar{\mu}})} ,\\
 \Vert  \bar{u} - (A\,\cdot +b) \Vert_{L^q(Q_{\bar{\mu}} \setminus F)} &\le \Vert \hat{u} - (A\,\cdot +b) \Vert_{L^q(Q_{\bar{\mu}})}\le  C \mu^{\frac{2}{q}-\frac{2}{p}+1}\Vert e(\hat{u}) \Vert_{L^p(Q_{\bar{\mu}})} 
 \end{align*}
 for $q \le \frac{2p}{2-p}$.  Then  the second part of \eqref{eq:modi2} holds for $p \in (1,2)$ by \eqref{eq:modi3} and the case $p=1$ directly follows.  Likewise, the first part of \eqref{eq:modi2} also holds for $q \in [1,\infty)$ since $p \in [1,2)$. \eop

\subsection{General case}

To prove the local estimate for a general function we consider a sequence of modifications and show that the properties in Theorem \ref{th: korn1} can be recovered in the limit.  

\noindent {\em Proof of Theorem \ref{th: korn1}.} Let $p\in (1,2)$, $q \in [1,\infty)$ and let $u \in SBD^2(Q_\mu)$ be given. By Theorem \ref{lemma: korn2} for $\delta= \frac{1}{n}$ we obtain  modifications $\bar{u}_n$ and exceptional sets $G_n$ such that \eqref{eq:modi} holds. Moreover, we find $A_n \in \R^{2 \times 2}_{\rm skew}$,  $b_n \in \R^2$ as well as exceptional sets $F_n \subset Q_\mu$ with ${\cal H}^1(\partial F_n) \le C{\cal H}^1(J_{{u}})$, $|F_n| \le C({\cal H}^1(J_{{u}}))^2$ such that
\begin{align}\label{eq: i,ii}
\begin{split}
(i) & \ \ \Vert \bar{u}_n - (A_n\,\cdot + b_n) \Vert_{L^q(Q_{\bar{\mu}} \setminus F_n)}\le C \mu^{\frac{2}{q}} \big(\Vert e(u) \Vert_{L^2( Q_{{\mu}})} + \tfrac{1}{n} \big),\\
(ii) & \ \ \Vert \nabla \bar{u}_n - A_n \Vert_{L^p(Q_{\bar{\mu}} \setminus F_n)}\le  C \mu^{\frac{2}{p}-1}\big(\Vert e(u) \Vert_{L^2( Q_{{\mu}})} + \tfrac{1}{n} \big), 
\end{split}\end{align}
where  $Q_{\bar{\mu}} = (-\bar{\mu},\bar{\mu})^2$ with $\bar{\mu} = \max\lbrace \mu - C{\cal H}^1(J_u) ,0 \rbrace$ independently of $n$.  Define $H_n = F_n \cup G_n$ and observe $|H_n| \le C({\cal H}^1(J_u))^2$, ${\cal H}^1(\partial H_n) \le C{\cal H}^1(J_u)$ for a sufficiently large constant. Then by Theorem \ref{clea-th: set per}  we find a set of finite perimeter $F \subset Q_\mu$ with ${\cal H}^1(\partial^* F) \le C{\cal H}^1(J_u)$, $|F| \le C({\cal H}^1(J_u))^2$ such that $\chi_{H_n} \to \chi_F$ in measure for $n \to \infty$ for a not relabeled subsequence. 

Moreover, letting $v_n := (\bar{u}_n - (A_n \,\cdot + b_n))\chi_{Q_{\bar{\mu}} \setminus H_n} \in GSBV^p(Q_\mu)$ and using \eqref{eq:modi}, \eqref{eq: i,ii} we apply Ambrosio's compactness result in $GSBV$ (see Theorem \ref{clea-th: compact}) to find a function $v \in GSBV^p(Q_\mu)$ such that passing to a further (not relabeled) subsequence we obtain $v_n \to v$ a.e. and  $\nabla v_n \rightharpoonup \nabla v$ weakly in $L^p$. In particular, we derive by Fatou's lemma
\begin{align*}
(i) & \ \ \Vert v \Vert_{L^q(Q_{\bar{\mu}} \setminus F)} \le \liminf\nolimits_{n \to \infty} \Vert v_n \Vert_{L^q(Q_\mu)} \le C \mu^{\frac{2}{q}} \Vert e(u) \Vert_{L^2( Q_\mu)},\\
(ii) & \ \ \Vert \nabla v \Vert_{L^p(Q_{\bar{\mu}} \setminus F)}\le \liminf\nolimits_{n \to \infty} \Vert \nabla v_n \Vert_{L^p(Q_\mu)} \le C  \mu^{\frac{2}{p}-1} \Vert e(u) \Vert_{L^2( Q_\mu)}.
\end{align*}
Consequently, to finish the proof it suffices to show $v = (u-a)\chi_{Q_{\bar{\mu}} \setminus F}$ for some infinitesimal rigid motion $a=a_{A,b}$. (Observe that as before the assertion then holds also for $p=1$.)

Possibly passing to a further subsequence we can assume $\chi_{H_n} \to \chi_{F}$ pointwise a.e. and thus we find a measurable set $B$ with $|B|>0$ such that $B \subset Q_{\bar{\mu}} \setminus F$, $B \subset  Q_{\bar{\mu}} \setminus H_n$ (up to a set of negligible measure) for $n$ large enough. By \eqref{eq:modi}, \eqref{eq: i,ii}(i) and H\"older's inequality this implies 
$$\Vert A_n \, \cdot +b_n \Vert_{L^1(B)} \le C \big(\Vert e(u) \Vert_{L^2( Q_{{\mu}})} + \tfrac{1}{n} \big)+ C\Vert \bar{u}_n \Vert_{L^1(B)} \le C$$ 
for $C=C(\mu)>0$ large enough. Consequently, we obtain $A_n \to A$, $b_n \to b$ for some $A \in \R^{2 \times 2}_{\rm skew}$ and $b \in \R^2$. As $\bar{u}_n \chi_{Q_{\bar{\mu}} \setminus H_n} \to u \chi_{Q_{\bar{\mu}} \setminus F}$ a.e. by \eqref{eq:modi} and $v_n \to v =v\chi_{Q_{\bar{\mu}} \setminus F} $ a.e., we conclude $v = (u - (A \, \cdot + b))\chi_{Q_{\bar{\mu}} \setminus F}$. 

As $GSBV^p(Q_\mu)$ is a vector space (see \cite[Proposition 2.3]{DalMaso-Francfort-Toader:2005}), we then get $u \chi_{Q_{\bar{\mu}} \setminus F} \in GSBV^p(Q_\mu)$. 
To see the additional statement that $u \chi_{Q_{\bar{\mu}} \setminus F} \in SBV^p(Q_\mu)$ if $u \in L^\infty(Q_\mu)$, we observe that $|A_n|,|b_n| \le C$ and $\Vert \bar{u}_n \Vert_{L^\infty(Q_{\bar{\mu}}\setminus G_n)} \le \Vert u \Vert_\infty$ (see Theorem \ref{rig-th: derive prop}) imply $\Vert v_n \Vert_\infty \le C$ independently of $n \in \N$ and the claim follows from Theorem \ref{clea-th: compact}.  \eop

 \section{Estimate at the boundary}\label{sec: bound}
 
In this section we give a refined estimate which holds up to the boundary of Lipschitz sets. This together with a standard covering argument will then lead to the proof of the main theorem. We first give an elementary estimate about the difference of infinitesimal rigid motions which we state in arbitrary space dimensions. 

\begin{lemma}\label{lemma: rigid}
Let $p \in [1,\infty)$ and $\bar{c}>0$. Then there is a constant $C=C(p,\bar{c})$ such that for all $x \in \R^d$, $R >0$ and measurable $\Omega \subset Q^x_R :=x+(-R,R)^d$ with $|\Omega|\ge \bar{c}R^d$ and all affine mappings $a: \R^d \to \R^d$ one has
$$\Vert a \Vert_{L^p(Q^x_R)}\le  C\Vert a \Vert_{L^p(\Omega)}.$$
\end{lemma}

Although similar estimates have already been used (see e.g. \cite{Chambolle-Conti-Francfort:2014, Friedrich:15-1}) we include the elementary proof here for the sake of completeness. 

\Proof We first note that by an elementary translation argument it suffices to consider cubes  $Q^0_R = (-R,R)^d$  centered at the origin. Assume the statement was false. Then there would be sequences $(R_k)_k$, $(\Omega_k)_k$ with $\Omega_k \subset Q^0_{R_k}$, $|\Omega_k| \ge \bar{c} R_k^d$ and a sequence of affine mappings $(a_k)_k$ with 
$$\Vert a_k \Vert_{L^p(Q^0_{R_k})}>  k\Vert a_k \Vert_{L^p(\Omega_k)}.$$
We define $c_k(x) = a_k(R_k x)$ as well as $\Omega'_k = \frac{1}{R_k} \Omega_k$ and obtain by transformation
$$\Vert c_k \Vert_{L^p(Q^0_{1})}>  k\Vert c_k \Vert_{L^p(\Omega'_k)}.$$ 
Then we define the affine mappings $d_k = \frac{c_k}{\Vert c_k \Vert_{L^p(Q^0_{1})}}$ and derive
$$1=\Vert d_k \Vert_{L^p(Q^0_{1})}>  k\Vert d_k \Vert_{L^p(\Omega'_k)}.$$
As $(d_k)_k$ are affine, we find that $\Vert d_k\Vert_{W^{1,p}(Q^0_1)}$ is uniformly bounded and thus a compactness result yields (after passage to a not relabeled subsequence) $d_k \to d$ in $L^p(Q^0_1)$ for some affine mapping $d$. Moreover, there is a measurable function $f$ with $f \ge 0$, $\Vert f \Vert_{L^1(Q^0_1)}\ge \bar{c}$ such that $\chi_{\Omega'_k} \rightharpoonup^* f$ weakly in $L^\infty(Q^0_1)$. Consequently, we find $1 = \Vert d \Vert_{L^p(Q^0_{1})}$ and $0 = \Vert d \cdot f\Vert_{L^1(Q^0_1)}$ which  gives a contradiction. \eop

We are now in a position to give the boundary estimate.  

 \begin{theorem}\label{th: boundary}
Let $p \in [1,2)$, $q \in [1,\infty)$.  Let $\mu>0$ and $\psi: (-2\mu,2\mu) \to [\mu,\infty)$ Lipschitz with $\Vert \psi'\Vert_\infty \le \bar{c}$ and $\inf \psi = \mu$. Let 
\begin{align}\label{eq: U}
 \begin{split}
 U &= \lbrace (x_1,x_2): - 2\mu < x_1 < 2\mu, \ -2\mu \le x_2 \le \psi(x_1) \rbrace, \\
U' &= \lbrace (x_1,x_2): - \mu < x_1 < \mu, \ -\mu \le x_2 \le \psi(x_1) \rbrace.
 \end{split}
\end{align}
Then there is a constant $C=C(p,q,\bar{c})$ independent of $\mu$ such that for all $u \in SBD^2(U)$ there is a set of finite perimeter $G \subset U$ with ${\cal H}^1(\partial^* G) \le C{\cal H}^1(J_u)$, $|G| \le C({\cal H}^1(J_u))^2$ and for suitable $A \in \R^{2 \times 2}_{\rm skew}$, $b \in \R^2$
\begin{align}\label{eq: main est***}
\begin{split} 
(i) & \ \ \Vert u - (A\,\cdot + b) \Vert_{L^q(U' \setminus G)}\le C \mu^{\frac{2}{q}} \Vert e(u) \Vert_{L^2(U)},\\
(ii) & \ \ \Vert \nabla u - A \Vert_{L^p(U' \setminus G)}\le C \mu^{\frac{2}{p}-1} \Vert e(u) \Vert_{L^2(U)}.
\end{split}
\end{align}

 \end{theorem}
 
\Proof Recall the definition of the sets ${\cal Q}_i$, $i \in \N$, and the enlarged squares  $Q   \subset Q'' \subset Q'$  in \eqref{eq: enlarged squares}, \eqref{eq:enlarged}. Moreover, by $d(B)$ we again denote the diameter of a set $B \subset \R^2$.   Let ${\cal Q}_W \subset \bigcup_{i \ge 1} {\cal Q}_i$ be a Whitney-type covering of $U$, i.e. $\bigcup_{Q \in {\cal Q}_W} Q' = U$  such that (cf. e.g. \cite{Acosta, Federer:1969, Stein:1970})
\begin{align}\label{eq: wit}
\begin{split}
 (i) & \ \ d(Q) \le \dist(Q,\partial U) \le Cd(Q) \ \ \text{ for all } \ Q \in {\cal Q}_W,\\
 (ii) & \ \ \# \lbrace Q \in {\cal Q}_W: x \in Q'\rbrace \le N \ \ \text{ for all } \ x \in U,\\
 (iii) & \ \ Q_1' \cap Q_2' \neq \emptyset \ \text{ for } \ Q_1,Q_2 \in {\cal Q}_W \  \Rightarrow \ \tfrac{1}{C} d(Q_1) \le d(Q_2) \le Cd(Q_1).
\end{split}
\end{align} 
Moreover, we consider a corresponding partition of unity $(\varphi_Q)_{Q \in {\cal Q}_W} \subset C^\infty(U)$ with $\sum_{Q \in {\cal Q}_W} \varphi_Q(x)  = 1$ for $x \in U$ and  
 \begin{align}\label{eq: prop3***}
\begin{split}
(i) & \  \ Q \subset {\rm supp}(\varphi_Q) \subset Q'' \ \text{ for all } Q \in {\cal Q}_{W},\\
(ii)& \ \ \Vert \nabla \varphi_Q\Vert_\infty \le C d(Q)^{-1} \ \text{ for all } \ Q \in {\cal Q}_{W}
\end{split}
\end{align}
for a universal constant $C>0$.  Let 
\begin{align}\label{eq: calcalB}
{\cal B} = \lbrace Q \in {\cal Q}_W: {\cal H}^1(Q' \cap J_u) \ge \hat{c}d(Q) \rbrace
\end{align}
be the `bad' squares for some $\hat{c}>0$ sufficiently small to be specified below. For each enlarged square $Q' =  p + (-r,r)^2$, $Q \in {\cal B}$, we define $P_Q = (p + (-r,r) \times (-r,\infty)) \cap U$. Employing \eqref{eq: wit}(i) and using $\Vert \psi' \Vert_\infty \le \bar{c}$ we then observe that  ${\cal H}^1(\partial P_Q) \le C d(Q)$ for $C=C(\bar{c})$ and thus ${\cal H}^1(\partial P_Q) \le C{\cal H}^1(Q' \cap J_u )$ for some $C=C(\bar{c},\hat{c})$. Letting $P = \bigcup_{Q \in {\cal B}} \overline{P_Q}$ we obtain by \eqref{eq: wit}(ii)
\begin{align}\label{eq: P}
{\cal H}^1(\partial P) \le CN{\cal H}^1(J_u)
\end{align}
and using the isoperimetric inequality we also find $|P| \le C({\cal H}^1(J_u))^2$. We let $V = U' \setminus P$.

Observe that we can assume $Q \notin {\cal B}$ for all $Q \in {\cal Q}_W$
with $Q \cap (- \mu,\mu) \times \lbrace 0 \rbrace \neq \emptyset$. In fact, these squares satisfy $d(Q)\ge c\mu$. Consequently, if $Q \in {\cal B}$, we find ${\cal H}^1(J_u) \ge c\hat{c}\mu$ and 
in this case the claim of the theorem holds with the choice $G = U$ if in the assertion $C$  is chosen large enough.

\vspace{-0.2cm}
\begin{figure}[H]
\centering
\begin{overpic}[width=0.65\linewidth,clip]{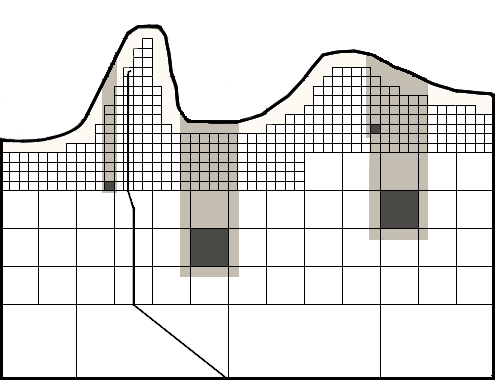}
\put(65,170){{$x$}}
\put(126,4){{$0$}}
\put(92,10){{$\gamma$}}
\end{overpic}
\caption{\small Illustration of a part of $U$ and ${\cal Q}_W \cap \bigcup^3_{i=1} {\cal Q}_i$. The squares in ${\cal B}$ are depicted in dark gray and the corresponding set $P$ in light gray. Moreover, a John curve $\gamma$ connecting $x$ with $0$ is sketched.} \label{bild}
\end{figure}

We now see that $V$ is a John domain with center $0$ and a constant only depending on $\bar{c}$. In fact, fix some $x = (x_1,x_2) \in V$ and $Q \in {\cal Q}_W$ such that
$x \in  \overline{Q}$.  We consider a vertical chain ${\cal C}_1 = \lbrace Q^1_1 = Q,\ldots, Q^1_{n_1}\rbrace$ of squares in ${\cal Q}_W$ intersecting
$\lbrace x_1 \rbrace \times [0,x_2]$ together with a horizontal chain ${\cal C}_2 = \lbrace Q^2_1 = Q^1_{n_1},\ldots, Q^2_{n_2}\rbrace$ of squares intersecting 
$[x_1,0] \times \lbrace 0 \rbrace$ such that $\overline{Q^j_{k}} \cap \overline{Q^j_{k+1}} \neq \emptyset$ for $1\le k \le n_j-1$, $j=1,2$.

Now in view of \eqref{eq: wit}, we see that $d(Q^j_{k_1}) \le d(Q^j_{k_2})$ for all $1 \le k_1 \le k_2 \le n_j$, $j=1,2$, and 
$d(Q^j_{k_1})\le \theta d(Q^j_{k_2})$ for all $k_2 \ge k_1 + l$ 
for some $l=l(\bar{c}) \in \N$. Consequently, it is elementary to construct
a curve $\gamma$ starting in $x$, ending in $0$ and intersecting the midpoints of the squares in ${\cal C}_1 \cup {\cal C}_2$ such that 
the condition given in Definition \ref{def: chain} holds (cf. Figure \ref{bild}).

Let ${\cal G} = {\cal Q}_W \setminus {\cal B}$. For each $Q \in {\cal G}$ we apply Theorem \ref{th: korn1} on $Q'$ to find  infinitesimal rigid motions $a_Q = a_{A_Q,b_Q}$ and exceptional sets $F_Q$ so that by \eqref{eq: main est} 
\begin{align}\label{eq: main est3}
\begin{split}
(i) & \ \ \Vert u - (A_Q\,\cdot + b_Q) \Vert_{L^q(Q'' \setminus F_Q)}\le C d(Q)^{\frac{2}{q}} \Vert e(u) \Vert_{L^2(Q')},\\
(ii) & \ \ \Vert \nabla u - A_Q \Vert_{L^p(Q'' \setminus F_Q)}\le C d(Q)^{\frac{2}{p}-1} \Vert e(u) \Vert_{L^2(Q')}.
\end{split}
\end{align}   
(For $\hat{c}$ sufficiently small in \eqref{eq: calcalB} we can in fact assume that $Q''$ is contained in the shrinked square   given by Theorem  \ref{th: korn1}.) Moreover, by \eqref{eq: R2} and \eqref{eq: wit}(ii) we get that  $F:= \bigcup_{Q \in {\cal G}} F_Q$ fulfills ${\cal H}^1(\partial^* F) \le C{\cal H}^1(J_u)$ and $|F| \le C({\cal H}^1(J_u))^2$.

We now estimate the difference of the infinitesimal rigid motions. Consider some $Q \in {\cal G}$ and let ${\cal N}(Q) =  \lbrace \hat{Q} \in {\cal G} \setminus \lbrace Q \rbrace: Q'' \cap \hat{Q}'' \neq \emptyset \rbrace$. Recall   $\frac{1}{C}d(\hat{Q}) \le d(Q) \le Cd(\hat{Q})$ for all $\hat{Q} \in {\cal N}(Q)$ by \eqref{eq: wit}(iii) which also implies $\# {\cal N}(Q) \le C$ for some $C>0$ large enough. 
Since the covering $\mathcal{Q}_W$ consists of dyadic squares, $Q'' \cap \hat{Q}''$ contains a ball $B$ with radius larger than $cd(Q)$ for some small $c>0$. In view of \eqref{eq: R2}, choosing $\hat{c}$ in \eqref{eq: calcalB} sufficiently small, we find   that $|\hat{F} \cap B|\le \frac{1}{2}|B|$, where $\hat{F} = F_{\hat{Q}} \cup F_Q$. Therefore, by  \eqref{eq: main est3}(i) for $p=q$, \eqref{eq: wit}(iii) and the triangle inequality we derive 
$$\Vert a_{Q} - a_{\hat{Q}}\Vert_{L^p(B \setminus \hat{F})}   = \Vert (A_{Q}\, \cdot + b_{Q}) - (A_{\hat{Q}}\, \cdot + b_{\hat{Q}})\Vert_{L^p(B \setminus \hat{F})} \le  C d(Q)^{\frac{2}{p}}  \Vert e(u) \Vert_{L^2(Q' \cup \hat{Q}')}$$
and thus by Lemma \ref{lemma: rigid}
\begin{align}\label{eq: diff1}
\Vert a_{Q} - a_{\hat{Q}}\Vert^p_{L^p(Q')} \le C d(Q)^{2}  \Vert e(u) \Vert^p_{L^2(Q' \cup \hat{Q}')}
\end{align} 
for some $C=C(p)$. Let $N_Q = \bigcup_{\hat{Q} \in {\cal N}(Q)}  \hat{Q}' \cup Q'$ and observe that by  \eqref{eq: wit} each $x \in U$ is contained in a bounded number of different sets $N_Q$. Moreover, we observe that $\sum_{Q\in {\cal G}}d(Q)^2 \le C|V| \le C\mu^2$.
Summing over all squares, recalling $\# {\cal N}(Q)\le C$ and using H\"older's inequality we then find
\begin{align}\label{eq: next}
 \begin{split}
\hspace{-0.2cm}&\sum\nolimits_{Q \in {\cal G}} \sum\nolimits_{\hat{Q} \in {\cal N}(Q)}   d(Q)^{-p}\Vert a_{Q} -a_{\hat{Q}}\Vert^p_{L^p(Q')} \le C\sum\nolimits_{Q \in {\cal G}} d(Q)^{2-p}  \Vert e(u) \Vert^p_{L^2(N_Q)}\\
& \ \ \le C\big(\sum\nolimits_{Q \in {\cal G}}d(Q)^2\Big)^{1-\frac{p}{2}}\big(\sum\nolimits_{Q \in {\cal G}} \Vert e(u) \Vert^2_{L^2(N_Q)}\Big)^{\frac{p}{2}} \le C\mu^{2-p}\Vert e(u) \Vert^p_{L^2(U)}.
\end{split}
\end{align}
We observe that $\sum\nolimits_{Q \in {\cal G}} \varphi_Q(x) = 1$ for all $x \in V$. In fact, we recall that $(\varphi_Q)_{Q \in {\cal Q}_W}$ is a partition of unity and ${\rm supp}(\varphi_Q)\subset Q'' \subset  U \setminus V$ for all $Q \in {\cal B} = {\cal Q}_W \setminus {\cal G}$ by construction. Similarly as in the proof of Theorem \ref{lemma: korn2} we define 
$$ \bar{u}(x) = \sum\nolimits_{Q \in {\cal G}} \varphi_Q(x) (A_Q \, x + b_Q) = \sum\nolimits_{Q \in {\cal G}} \varphi_Q(x) a_Q(x)$$
for all $x \in V$. Clearly, $\bar{u}$ is smooth in $V$. Using $\sum\nolimits_{Q \in {\cal G}}  \nabla \varphi_{Q} = 0$ we find that 
\begin{align}\label{eq: deriva}
\nabla \bar{u} = \sum\nolimits_{\hat{Q} \in {\cal G}} \big( (a_{\hat{Q}}  -f) \otimes \nabla \varphi_{\hat{Q}} + \varphi_{\hat{Q}} A_{\hat{Q}}\big)
\end{align}
for all functions $f$. Consequently, letting $f(x) = a_Q(x)$ for $x \in Q \cap V$, $Q \in {\cal G}$, we derive using once more \eqref{eq: wit}(ii),(iii) and applying \eqref{eq: prop3***}, \eqref{eq: next}
\begin{align*}
\Vert e(\bar{u})\Vert^p_{L^p(V)}  &\le C\sum\nolimits_{Q \in {\cal G}}  d(Q)^{-p} \sum\nolimits_{\hat{Q} \in {\cal N}(Q)}   \Vert a_{Q} -a_{\hat{Q}}\Vert^p_{L^p(Q)} \le C \mu^{2-p}\Vert e(u) \Vert^p_{L^2(U)}.
\end{align*}
We compute for $q \ge 2$ using \eqref{eq: wit}(ii), \eqref{eq: prop3***}(i) and \eqref{eq: main est3}(i)
\begin{align}\label{eq: 1}
\Vert \bar{u} - u\Vert^q_{L^q(V\setminus F)} & \le C\sum\nolimits_{Q \in {\cal G}} \Vert a_Q - u\Vert^q_{L^q(Q''\setminus F)} \le C\sum\nolimits_{Q \in {\cal G}} d(Q)^{2}\Vert e(u)\Vert^q_{L^2(Q')} \notag\\& \le C\mu^{2}\Vert e(u)\Vert^q_{L^2(U)}. 
\end{align}
(The case $1 \le q <2$ follows similarly by H\"older's inequality.) Likewise, by \eqref{eq: deriva} for $f = u$ and \eqref{eq: main est3} for $q=p$ we find repeating the H\"older-type estimate in \eqref{eq: next}
\begin{align}\label{eq: 2}
\begin{split}
\Vert \nabla \bar{u} - \nabla u\Vert^p_{L^p(V\setminus F)} & \le C\sum\nolimits_{Q \in {\cal G}} d(Q)^{-p} \Vert a_{Q} - u\Vert^p_{L^p(Q''\setminus F)} \\ & \ \ \ +C\sum\nolimits_{Q \in {\cal G}}   \Vert \nabla u - A_{Q}\Vert^p_{L^p(Q''\setminus F)}   \\
&\le  C\sum\nolimits_{Q \in {\cal G}} d(Q)^{2-p}\Vert e(u)\Vert^p_{L^2(Q')} \le C\mu^{2-p}\Vert e(u)\Vert^p_{L^2(U)}.
\end{split}
\end{align}
As $\bar{u}$ is smooth in $V$ and $V$ is a John domain with constant only depending on $\bar{c}$, we can apply Theorem \ref{th: kornsobo} and we find $A \in \R^{2 \times 2}_{\rm skew}$, $b\in\R^2$ such that by a scaling argument 
$$\mu^{-\frac{2}{q}-1 + \frac{2}{p}}\Vert \bar{u} - (A\,\cdot + b) \Vert_{L^q(V)} + \Vert \nabla \bar{u} - A \Vert_{L^p(V)}\le C  \Vert e(\bar{u}) \Vert_{L^p(V)} \le C\mu^{\frac{2}{p}-1} \Vert e(u) \Vert_{L^2(U)}$$
for $C=C(p,q,\bar{c})$.  We now define $G = F\cup P$ and by \eqref{eq: P} and the remark below \eqref{eq: main est3} we obtain $|G| \le C({\cal H}^1(J_u))^2$ as well as ${\cal H}^1(\partial^* G) \le C{\cal H}^1(J_u)$.  Finally, \eqref{eq: main est***} follows
from \eqref{eq: 1} and \eqref{eq: 2}. \eop 

\begin{rem}
{\normalfont

Similarly as in the local estimate considered in Section \ref{sec: local} one can show that the displacement field restricted to $U' \setminus G$ is an element of $GSBV^p$ or $SBV^p$, respectively. As this property will not be needed in the following, we have omitted the proof. 

}
\end{rem}

 \section{Proof of the main result and application}\label{sec: main}
 
 \subsection{Proof of Theorem \ref{th: mainkorn}}
We now combine the local estimate in Theorem \ref{th: korn1}, the boundary estimate (Theorem \ref{th: boundary})  and a standard covering argument to prove the main result. A similar argument may be found, e.g. in \cite{Chambolle-Conti-Francfort:2014}, where an inequality of Korn-Poincar\'e type is derived.

\noindent {\em Proof of Theorem \ref{th: mainkorn}}. We first choose finitely many $U_1,\ldots, U_n$ being of the form given in \eqref{eq: U} (possibly after application of an affine isometry)
such that $\partial \Omega$ is covered by $U_1',\ldots,U_n'$. Moreover, we cover $\Omega \setminus \bigcup^n_{i=1} U_i'$ with squares $U_{n+1}',\ldots,U_{m}'$ such that the squares $U_{n+1} := 2U_{n+1}',\ldots, U_m:=2U_{m}'$
of double size are still contained in $\Omega$.

By a similar reasoning as in the proof of Theorem \ref{th: boundary} we may suppose that  ${\cal H}^1(J_u) \le \hat{c}$ for some $\hat{c} = \hat{c}(p,q,\Omega)$ to be specified below as otherwise we can choose $F = \Omega$ in Theorem \ref{th: mainkorn}. 
We now apply Theorem \ref{th: korn1} and Theorem \ref{th: boundary}, respectively, on the sets $(U_i)^m_{i=1}$ and obtain infinitesimal rigid motions $a_i = a_{A_i,b_i}$ as
well as exceptional sets $F_i \subset U_i$ such that
\begin{align}\label{eq: last}
\sum^m_{i=1} \big(\Vert u - a_i\Vert_{L^q(U_i'\setminus F_i)}  + \Vert \nabla u - A_i\Vert_{L^p(U_i'\setminus F_i)}\big)
\le C\sum^m_{i=1} \Vert e(u)\Vert_{L^2(U_i)} \le C\Vert e(u)\Vert_{L^2(\Omega)}
\end{align}
for some $C=C(p,q,\Omega)$. In fact, selecting $\hat{c}$ sufficiently small we get that the shrinked squares given in Theorem \ref{th: korn1} contain $U_i'$ for $i=n+1,\ldots,m$ (cf. \eqref{eq: main est3} for a similar argument).

Define $F = \bigcup^m_{i=1} F_i$ and observe that $|F| \le C({\cal H}^1(J_u))^2$  as well as ${\cal H}^1(\partial^* F) \le C{\cal H}^1(J_u)$ 
follow from \eqref{eq: R2} and the similar estimate for the sets at the boundary (see before \eqref{eq: main est***}). Moreover, we can choose $\hat{c}$ so small such that $|F| \le \frac{1}{2}\eta$, where 
$$\eta:= \min\lbrace |U_i' \cap U_j'|:  U_i', U_j', i \neq j, \text{ with } U_i' \cap U_j' \neq \emptyset\rbrace.$$
Obviously, $\eta$ only depends on $\Omega$. Consequently, we obtain $|(U_i' \cap U_j') \setminus F| \ge \frac{1}{2}\eta$ for all $U_i'$, $U_j'$, $i \neq j$, with $U_i' \cap U_j' \neq \emptyset$. As $\Omega$ is connected, we then find by Lemma \ref{lemma: rigid} and \eqref{eq: last}
$$\max_{1 \le i,j \le m} \big(\Vert a_i - a_j\Vert_{L^q(\Omega)} + \Vert A_i - A_j \Vert_{L^p(\Omega)}\big) \le C\Vert e(u) \Vert_{L^2(\Omega)}$$
for a constant depending only on $p,q$, $\eta$ and $m$. Recalling \eqref{eq: last} and the fact that $\eta$, $m$ only depend on $\Omega$, we finally obtain \eqref{eq: main estmain} for, e.g., $A = A_1$ and $b=b_1$. \eop

\subsection{Relation between $SBV$ and $SBD$ functions}
Finally, we present a consequence of our main result concerning the relation between $SBV$ and $SBD$ functions. We briefly recall that the typical examples
for functions lying in $BD$ but not in $BV$ or likewise lying in $SBD^p$ but not in $SBV^p$, $p>1$, are based on the idea to cut out small balls  and to choose the displacement field appropriately on these sets (see e.g.  \cite{Ambrosio-Coscia-Dal Maso:1997, Conti-Iurlano:15}). 
The following result shows that this construction essentially describes the only way to obtain functions of bounded deformation which do not have bounded variation. In particular, we see that for each function in $SBD^2 \cap L^\infty$ there is a modification $\tilde{u}$ in $SBV$ such that $\lbrace u \neq \tilde{u} \rbrace$  is an arbitrarily small set.  

\begin{theorem}\label{th: appl}
 Let $\eps > 0$ and let $\Omega \subset \R^2$ open, bounded with Lipschitz boundary. Then for every $u \in SBD^2(\Omega)$ we find an exceptional
 set $F$ with $|F| \le \eps$ and ${\cal H}^1(\partial^* F)< + \infty$ such that $u \chi_{\Omega\setminus F} \in GSBV^p(\Omega)$ for all $p<2$. If in addition $u \in L^\infty(\Omega)$, we obtain $u \chi_{\Omega\setminus F} \in SBV^p(\Omega)$.
\end{theorem}

\Proof The statement follows from Theorem \ref{th: korn1} by an additional covering argument. Assume first $u \in SBD^2(\Omega) \cap L^\infty(\Omega)$. Recalling \eqref{eq: enlarged squares} we cover $\Omega$ with squares in ${\cal Q}_i$ for $s_i \ll \eps$ to be specified below. We define the bad squares ${\cal B} = \lbrace Q \in {\cal Q}_i: Q' \not\subset \Omega \ \text{ or } \   {\cal H}^1(J_u \cap Q') \ge \bar{c}s_i \rbrace$ for a constant $\bar{c}>0$  to be specified below.  We let 
$G = \Omega\cap \bigcup_{Q \in {\cal B}} \overline{Q}$ and observe ${\cal H}^1(\partial G) \le C{\cal H}^1(J_u) + C{\cal H}^1(\partial \Omega) < + \infty$ as well as $|G| \le Cs_i({\cal H}^1(J_u) + {\cal H}^1(\partial \Omega))$. 

Choosing $\bar{c}$ 
sufficiently small we can apply Theorem \ref{th: korn1} on each enlarged square $Q', Q \in {\cal Q}_i \setminus {\cal B}$, and obtain exceptional sets $F_Q$ such that $\big( u \chi_{Q'' \setminus F_Q}\big)|_{Q'}\in SBV^p(Q')$ for all  $Q \in {\cal Q}_i \setminus {\cal B}$. (In fact, for $\bar{c}$ small we can assume that the shrinked square  given in Theorem \ref{th: korn1} contains $Q''$.)

Letting $F = \bigcup_{Q} F_Q \cup G$ we find by \eqref{eq: R2} that
${\cal H}^1(\partial^* F) < +\infty$ and $|F| \le Cs_i$ for a constant depending only on $\Omega$ and $u$. Thus, $|F| \le \eps$ for $i \in \N$ sufficiently large. Then defining $\bar{u}= u \chi_{\Omega \setminus F}$ in $\Omega \setminus G$ we derive $\bar{u} \in SBV^p(\Omega \setminus G)$ for all $p \in [1,2)$. Observe that
\begin{align}\label{eq: trace-exp}
D(u \chi_{\Omega \setminus F}) = D\bar{u} + (\bar{u} \otimes \xi_G) \, {\cal H}^1|_{\partial G \cap \Omega}
\end{align}
in $\Omega$, where $\xi_G$ denotes the inner normal of $G$ (see e.g. \cite[Theorem 3.87]{Ambrosio-Fusco-Pallara:2000}). As $\Vert u \Vert_\infty < + \infty$, this implies $u \chi_{\Omega \setminus F} \in SBV^p(\Omega)$. Likewise, in the general case we consider $\phi \in C^1(\R^2)$ with the support of $\nabla \phi$ compact and find $\big(\phi(u \chi_{\Omega \setminus F}) \big)|_{\Omega \setminus G} \in SBV^p(\Omega \setminus G)$. Then we repeat the argument in \eqref{eq: trace-exp} to conclude $\phi(u \chi_{\Omega \setminus F}) \in SBV^p(\Omega)$. \eop

The above result can also be interpreted as an approximation result for $SBD$ functions. On the one hand, it is weaker than standard density results, see e.g. \cite{Chambolle:2004}, as 
it does not lead to a fine estimate for the surface energy. On the other hand, whereas  in results based on interpolation arguments the approximating sequences typically only converge in $L^p$,
in the present context we see that the functions already coincide up to a set of arbitrarily small measure.

\bigskip
 
\noindent \textbf{Acknowledgements} This work has been funded by the Vienna Science and Technology Fund (WWTF)
through Project MA14-009.


 \typeout{References}

\end{document}